\documentclass{amsart}
\usepackage{amssymb,amscd,verbatim}
\usepackage{enumerate}
\usepackage{graphicx}
\newcommand\datver[1]{\def\datverp
{\par\boxed{\boxed{\text{Version: #1; Run: \today}}}}}
\datver{1.1; Revised: 01/16/2004}

\newcommand{\VV}{\mathcal V}
\newcommand{\pa}{\partial}
\newcommand{\CI}{\mathcal C^{\infty}}
\newcommand\CIc{\mathcal{C}^\infty_c}
\newcommand{\Diff}{\operatorname{Diff}}

\newcommand{\rx}{{x}}

\newcommand\PP{\mathbb P}
\newcommand\RR{\mathbb R}

\newcommand{\maG}{\mathcal G}
\newcommand{\maH}{\mathcal H}

\newcommand{\maL}{\mathcal L}

\newcommand{\mfA}{\mathfrak A}
\newcommand{\mfL}{\mathfrak L}

\newcommand{\ie}{{\em i.e., }}

\newtheorem{theorem}{Theorem}[section]
\newtheorem{proposition}[theorem]{Proposition}
\newtheorem{corollary}[theorem]{Corollary}

\theoremstyle{definition}
\newtheorem{definition}[theorem]{Definition}
\theoremstyle{remark}
\newtheorem{remark}[theorem]{Remark}
\newtheorem{example}[theorem]{Example}

\numberwithin{figure}{section} \numberwithin{table}{section}

\newcommand\PSM[1]{\Psi_{\VV}^{#1}(M)}
\newcommand\mate[1]{{$ #1 $}}

\author[V. Nistor]{Victor Nistor}
\email{nistor@math.psu.edu}
\address{Pennsylvania State University,
Math. Dept., University Park, PA 16802}

\thanks{The author was supported in part by the NSF grants DMS
0209497 and DMS 991981.}

\begin{document}

\dedicatory\datverp

\title[Polyhedral domains]{Singular integral operators
on non-compact manifolds and analysis on polyhedral domains}

\begin{abstract}
We review the definition of a Lie manifold $(M, \VV)$ and the
construction of the algebra $\PSM{\infty}$ of pseudodifferential
operators on a Lie manifold $(M, \VV)$. We give some concrete
Fredholmness conditions for pseudodifferential operators in
$\PSM{\infty}$ for a large class of Lie manifolds $(M, \VV)$.
These Fredholm conditions have applications to boundary value
problems on polyhedral domains and to non-linear PDEs on
non-compact manifolds. As an application, we determine the
spectrum of the Dirac operator on a manifold with
multi-cylindrical ends.
\end{abstract}

\maketitle \tableofcontents

\section*{Introduction}

Partial diffferential equations on non-compact manifolds are a common
occurrence in Geometry, Group Representations, Mathematical Physics,
and other areas of Mathematics and Science. For example, conformally
compact manifolds and asymptotically flat manifolds were recently
considered in Quantum Gravity and in the study of the AdS-CFT
correspondence \cite{anderson:p03, Ashtekar, MazzeoQG, GraZwor, Lee,
SchroheQG}. One of the main technical issues of the analysis on
non-compact manifold $M_0$ is that an elliptic operator $P$ of order
$m$ with elliptic principal symbol is not necessarily Fredholm as an
operator $P : H^{s}(M_0) \to H^{s-m}(M_0)$. In particular, the
spectrum of such a $P$ needs not to be discrete.

Analysis on non-compact manifolds plays a role also in the analysis on
singular spaces, the non-compact space being the set of regular points
endowed with a suitable metric. An important class of singular spaces
is provided by polyhedral domains.

Analysis on polyhedral domains has many features that are not present
in the analysis on smooth domains. Several of these issues were
discussed for Lipschitz domains in \cite{JerisonKenig, MitreaTaylor,
Necas, Verchota}. However, polyhedral domains are not always Lipschitz
(recall the ``two-brick'' example \cite{VerVog1}). Moreover,
polyhedral domains are amenable to a more detailed analysis
\cite{Costabel, Dauge, Grisvard1, Grisvard2}. So far, however, this
more detailed analysis was devoted mostly to the case of polygonal
domains, and, occasionally, to the case of polyhedral domains in
space. See however the recent work of Verchota and Vogel on higher
dimensional polyhedra \cite{VerVog1, VerVog2}.

In this paper, we discuss the relevance of singular integral operators
in the analysis on non-compact manifolds and in the analysis on
polyhedral domains. This paper is largely based on joint results
with:\ Bernd Ammann, Constantin Bacuta, Robert Lauter, Alexandru
Ionescu, Marius Mitrea, Bertrand Monthubert, Andras Vasy, Alan
Weinstein, Ping Xu, and Ludmil Zikatanov \cite{AIN, aln1, aln2, alnv,
BNZ1, BNZ2, LMN1, LN1, MitreaNistor, NWX}. A central role in the above
papers is played by the concept of Lie manifold \cite{aln1} (their
definition is recalled in Definition \ref{def.LieM}) and by the
natural pseudodifferential operators acting on a Lie manifold. In
\cite{aln1}, the term ``manifold with a Lie structure at infinity''
was used instead of the term ``Lie manifold.''

We begin by recalling some results on boundary value problems that
motivate our interest in non-compact manifolds. Then we recall the
definition of a Lie manifold $(M, \VV)$, where $\VV$ is a suitable Lie
algebra of vector fields on $M$ and the construction of the Melrose
quantization $\PSM{\infty}$ of the algebra of differential operators
naturally associated to a Lie manifold. In addition to reminding some
of the necessary results from the above papers, we also prove some new
results. We introduce a special class of Lie manifolds (``type I Lie
manifolds'') in Section \ref{sec.FC} and we give some explicit
conditions for an operator $P \in \PSM{\infty}$ to be Fredholm if $(M,
\VV)$ is a Lie manifold.  Section \ref{sec.Ex} contains several
concrete examples of Lie manifolds and of the use of the Fredholmness
conditions for type I Lie manifolds. As a last application, we also
determine in Section \ref{sec.Sp} the essential spectrum of the Dirac
operator on a manifold with multi-cylindrical ends. For comparison,
let us recall that the essential spectrum of the Laplace operator on a
manifold with multi-cylindrical ends is $[0, \infty)$ \cite{LN1}
(solving a conjecture from \cite{MelroseScattering}).

The structure of this paper reflects, to a large extent, the structure
of my talk given at the ``Conference on Spectal Geometry of Manifolds
with Boundary and Decomposition of Manifolds,'' organized by
B. Booss-Bavnbek, G. Grubb, and K. Wojciechowski, whom I thank for
their efforts and for the opportunity to present my results. This
paper, however, contains more precise statements and several new
results. I also thank Bernd Ammann, Constantin B\v{a}cu\c{t}a, Craig
Evans, Alexandru Ionescu, Robert Lauter, and Irina Mitrea for useful
discussions. We shall write ``$:=$'' for ``the left hand side is equal
by definition to the right hand side.''

%%
%% TO INCLUDE   The paper is organized as follows:
%%

\section{Boundary value problems}

Let \mate{\Omega \subset \RR^n} be a {\em bounded, open} set with
boundary ${ \pa \Omega := \overline{\Omega} \smallsetminus
\Omega\,.}$. Let us consider on $\Omega$ ``simplest'' boundary
value problem, the Poisson problem
\begin{equation}\label{eq.BVP}
    \begin{cases} \; \Delta u = f&\\
    u\vert_{\pa \Omega} = g. &
\end{cases}
\end{equation}

A well known, classical results \cite{Evans, Taylor1} is the
following ``shift theorem'' (or ``regularity theorem'').

\begin{theorem}[Classical]\label{thm.reg0}\
If \mate{\pa \Omega} is {smooth}, then \mate{\tilde{\Delta}(u) =
(\Delta u, u \vert_{\pa \Omega})} defines an isomorphism
\begin{equation*}
     \tilde{\Delta} : H^{m + 2}(\Omega)  \to
     H^{m}(\Omega) \oplus H^{m + 3/2}(\pa \Omega),
\end{equation*}
for $m \in \RR$, $m \ge -1$.
\end{theorem}

Although this is not needed in this paper, let us notice, for the
interested reader, that the range of $m$ in the above theorem can be
improved to $m \in \RR$.  See for example \cite{BabuNistor1} and the
references therein. This improvement is relevant in some problems
arising in the applications of elasticity to Engineering
\cite{KoitSter}.

It follows right away from the above theorem that if $f$, $g$, and
\mate{\pa \Omega} are smooth, then \mate{u} is also smooth
(including the boundary). This is, however, not true in general if
$\pa \Omega$ is not smooth. In particular, the above theorem is
not true if \mate{\pa \Omega} is not smooth. Indeed, let us take
\mate{\Omega} to be the square \mate{(0, 1)^2} and $g = 0$ and let
us assume that $u$ is smooth. Then $\pa_x^2u(0, 0) = 0 =
\pa_y^2u(0, 0)$, and hence \mate{f(0,0) = \Delta u(0, 0) = 0}.
Thus no solution of our problem \eqref{eq.BVP} on the unit square
$\Omega = (0, 1)^2$ is smooth if $g = 0$ and $f(0, 0) \neq 0$. The
same problem arises on any polygonal domain. A detailed and far
reaching analysis of the above issues can be found in the
fundamental paper of Jerison and Kenig  \cite{JerisonKenig}, which
shows the exact range of applicability of the above theorem on a
Lipschitz domain in the plane.

From a practical point of view, the fact that the above ``shift
theorem'' does not extend directly to non-smooth domains is quite
inconvenient for applications. More precisely, if one want to solve an
elliptic partial differential equation using the finite element method
and a quasi-uniform mesh, the rate of convergence of the method is
governed by the smoothness of the solution. In particular, one
achieves only low orders of convergence using quasi-uniform meshes on
a polygon \cite{Wahlbin84}. The problem, however, is due to the use of
quasi-uniform meshes (and the use of the usual, isotropic Sobolev
spaces). Indeed, it was shown by Babu\v{s}ka already in the '70
\cite{Babuska70} that one can achieve the same rate of convergence as
for smooth domains, provided that one chooses correctly the finite
element space. See also \cite{BNZ1} and \cite{Raugel}.

In this paper we look at these issues from the point of view of Lie
manifolds. (We shall recall the definition of Lie manifolds and their
relevance to boundary value problems below.) Let us begin by
discussing the relatively simpler example of a polygonal domain (or,
more generally, of a domain whose boundary has conical points).

One of the most successful approaches so far is to use polar
coordinates \mate{(r,\theta)} around the vertices of a polygon.  For a
general domain with conical points, one uses generalized polar
coordinates. In the mathematical community this approach was pioneered
by Kondratiev \cite{Kondratiev67}, but in the Engineering community
the use of polar coordinates and of the Mellin transform is apparently
much older.

To explain this approach, let us consider the open angle \mate{\Omega
= \{ \theta \in (0, \alpha) \} } and polar coordinates, then
\begin{equation*}
    {\Delta = r^{-2}\big ( (r\pa_r)^2 + \pa_{\theta}^2 \big)}.
\end{equation*}
This suggests to look at differential operators on $\Omega$ of the
form
\begin{equation}\label{eq.delta}
    {\sum_{i+j \le m} a_{ij}(r, \theta) ({r}\pa_r)^i \pa_{\theta}^j}
\end{equation}
with \mate{a_{ij}} smooth. Operators of this type are called {\em
totally characteristic operators}, and they are defined on any
manifold with boundary (see below). The relevant part of the boundary
in our example is given by $r = 0$.

Analogously, we define then the {\em totally characteristic vector
fields} on $[0, \infty) \times (0, \alpha) \ni (r, \theta)$ to be the
vector fields $X$ of the form
\begin{equation}\label{eq.tot.char}
    X = a(r, \theta) {r}\pa_r + b(r, \theta) \pa_{\theta},
\end{equation}
where $a$ and $b$ are smooth functions. An important observation
\cite{meicm, MelroseScattering}, is that the totally characteristic
vector fields form a Lie algebra. This observation extends to
polygonal domains, domains with conical points, and, more generally,
to polyhedral domains.

Before explaining our use of vector fields, let us take a look at two
more examples. Let us consider the ``edge'' $\Omega \times \RR$, where
$\Omega = \{ \theta \in (0, \alpha) \}$, as above. If we consider
cylindrical coordinates \mate{(r, \theta, z)} in \mate{\RR^3}, then
the Laplace operator becomes
\begin{equation*}
    \Delta = r^{-2}\big ( (r\pa_r)^2 + \pa_{\theta}^2 +
    r^2\pa_z^2\big).
\end{equation*}
If we ignore the coefficient $r^{-2}$, we are lead to consider
differential operators generated by products of the derivatives $r
\pa_r$, $\pa_\theta$, and $r\pa_z$ (and smooth coefficients). The
differential operators of this kind that are vector fields are of the
form
\begin{equation}\label{eq.tot.char2}
    X = a(r, \theta, z) r\pa_r + b(r, \theta, z) \pa_{\theta} + c(r,
    \theta, z)r\pa_z,
\end{equation}
with $a$, $b$, and $c$ smooth functions. These vector fields
(``edge--type vector fields'') form also a Lie algebra.

\section{Lie algebras of vector fields }

The examples of the previous section, among others, have led Melrose
to formulate a program to study the analysis of differential operators
generated by suitable Lie algebras of vector fields \cite{meicm,
MelroseScattering}. Many important results in this program were
obtained by \cite{emm91, Lauter, LauterMoro, Mazzeo,
MelroseScattering, MelroseMendoza, ScSC, Schulze98, jaredduke}. In
this paper, however, we shall be mainly concerned with the approach to
this program developed in \cite{aln1, aln2, LMN1, LN1, Monthubert,
NWX}.

We shall consider a compact manifold with corners \mate{M}
together with a subspace \mate{\VV \subset \Gamma(TM)}, consisting
of vector fields {tangent} to all faces of \mate{M} and satisfying
certain conditions that make $(M, \VV)$ a ``Lie manifold.'' We
shall denote by $M_0$ the interior of $M$ and by $\pa M$ the set
of boundary points of $M$. In particular, $M_0 = M \smallsetminus
\pa M$. The following definition is essentially from
\cite{MelroseScattering}, but it was formalized in \cite{aln1}.

\begin{definition}\label{def.LieM}\
Let $M$ be a manifold with corners. A {\em Lie manifold} is a pair
$(M, \VV)$, where $\VV$ is a set of  vector fields tangent to all
faces of \mate{M} satisfying the following conditions:
\begin{enumerate}[(i)]
\item\ \mate{\VV} is closed under the Lie bracket \mate{[\;,\;]};
\item\ \mate{\CI(M)\VV = \VV};
\item $\VV$ is linearly generated locally in the neighborhood of each
point \mate{p \in M} by \mate{n} linearly independent vector fields
\mate{X_1, \ldots, X_n} with \mate{\CI(M)} coefficients.
\item\ If in the conditions above \mate{p \in M_0}, then the vector
fields \mate{X_1, \ldots, X_n}, locally generating \mate{\VV} around
\mate{p}, also give a local basis of \mate{T M} around $p$.
\end{enumerate}
\end{definition}

Condition (iii) means the following. For each $p\in M$ there exists an
open neighborhood $U$ of $p$ in $M$ and vector fields
$X_1,\ldots,X_n\in \VV$ such that for any $X \in \VV$, there exist
uniquely determined smooth functions $a_1, \ldots, a_n$ such that
\begin{equation}
	X = \sum a_i X_i \text{ on } U.
\end{equation}

The following remark is slightly less elementary, but it will be
useful in several places. It also explains the above definition.

\begin{remark}\label{rem.Algebroid}\
It follows from the last axiom that the integer $n$ appearing
above must be the same as the dimension of $M$. In particular,
$\VV$ is a $\CI(M)$--module isomorphic to a direct summand of the
free $\CI(M)$--module $\CI(M)^N$, for some $N$. That is $\VV$, is
a projective $\CI(M)$--module. Then the Serre-Swan theorem states
that there exists a vector bundle $A \to M$, unique up to
isomorphism, such that $\VV$ is isomorphic, as a $\CI(M)$-module
to $\Gamma(A)$, the space of sections of $A$.  The definition of
$\VV$ as a space of vector fields on $M$ and the naturality of
$A$, show that there exists a vector bundle map $\varrho : A \to
TM$, called {\em anchor map}, that endows $A$ with the structure of a
Lie algebroid. We shall therefore call $A$ the {\em Lie algebroid
associated to the Lie manifold $(M, \VV)$}.  Condition (iv) of
Definition \ref{def.LieM} is then equivalent to saying that
$\varrho$ is an isomorphism on the interior of $M$. See
\cite{aln1, LN1} for more details.
\end{remark}

In \cite{aln1}, the manifolds introduced in the above definition
were called ``manifolds with a Lie structure at infinity.''

Define \mate{\Diff_{\VV}(M)} to be the algebra of differential
operators on $M$ generated by $\VV$ and $\CI(M)$. The differential
operators in \mate{\Diff_{\VV}(M)} are the singular differential
operators we plan to study in this paper, due to their
applications to analysis on singular domains and on non-compact
manifolds.

Even if one is primarily interested in differential operators, in
order to invert them, one has to consider also integral kernel
operators. In our case, these integral kernel operators will be
pseudodifferential operators. To see their relevance for boundary
value problems, in particular, let us quickly recall the method of
layer potentials.

Let \mate{\Omega \subset \RR^n} be a bounded Lipschitz domain (for
example, a domain with piecewise $C^1$-boundary). Let $c_n^{-1}=
\omega_n (2-n)$, where $\omega_n$ is the surface of the unit sphere in
$\RR^n$), $\nu(y)$ is the outer unit normal, and $d\sigma(y)$ is the
induced measure on $\pa \Omega$. Then the operator
%\begin{equation*}
%    {
%    Sf(x) = \frac{c_n}{2-n} \int_{\pa \Omega} \frac{f(y)}{|y -
%    x|^{n-2}} d\sigma(y)}
%\end{equation*}
%and
\begin{equation*}
    Kf(x) = c_n \int_{\pa \Omega} \frac{(y-x) \cdot \nu(y)}
    {|y - x|^n} f(y) d\sigma(y),
\end{equation*}
can be used to determine the boundary value of the double layer
potential operator (the double layer potential operator is the
extension of the formula for $K$ to $x$ in the interior of $\Omega$,
while $y \in \pa \Omega$). If one can establish the invertibility of
\mate{\frac{1}{2}I + K} as a pseudodifferential operator on $\pa
\Omega$, then one obtains that the boundary value problem
\eqref{eq.BVP} has a solution for $f \in H^{1/2}(\pa \Omega)$, $g =
0$, which can be then used to obtain a solution for more general data.

If \mate{\pa \Omega} is smooth, then \mate{K} is a
pseudodifferential operator of order \mate{-1}. Hence
\mate{\frac{1}{2}I + K} is a Fredholm operator of index zero,
because \mate{\pa \Omega} is also compact. Therefore
\mate{\frac{1}{2}I + K} is invertible if, and only if, it is
injective or surjective. The injectivity can usually be checked
using energy methods. This completes our very brief summary of the
method of layer potentials for smooth domains.

The above reasoning does not extend directly to the case when
\mate{\pa \Omega} is not smooth, because $K$ may fail to be compact
\cite{FJL1, FJL2}. See also \cite{Verchota}. Nevertheless, for the
case of a polygon, $K$ is in a class of operators that is well
understood (the class of Hardy-type operators), see \cite{FJL2,
LewisParenti}. An approach to the study of Hardy-type operators is
provided by the operators in the ``$b$-calculus'' on \mate{\pa \Omega}
\cite{LewisParenti, MelroseScattering} to which the Hardy-type
operators are closely related. The $b$-calculus is the
pseudodifferential analog of totally-characteristic differential
operators. Using an iterative argument, one can show that the method
of layer potentials extends domains with conical points (hence to to
curvilinear polygons as well) \cite{FJL1, FJL2, LewisParenti,
MitreaNistor}. Let $r$ denote the distance the set of singularities on
the boundary (\ie the distance to the vertices, if our domain is a
polygon, or the distance to the conical points, if our domain is a
domain with conical points). Then define
\begin{equation*}
    r^{a}H^m_b(\Omega) := \{ u \in L^2_{loc}(\Omega),
    r^{-a-1+ |\alpha|}\pa^{\alpha} u \in L^2(\Omega),
    |\alpha| \leq m\}.
\end{equation*}

\begin{theorem}\label{thm.reg1}\
Let \mate{\Omega} be a polygon or a domain with conical points.
Then there exists $\eta > 0$ such that the map
\mate{\tilde{\Delta} (u) = (\Delta u, u\vert_{\pa \Omega})}
establishes an isomorphism
\begin{equation*}
    \tilde{\Delta} : r^aH^{m+2}_b(\Omega) \to
    r^{a-2}H^{m}_b(\Omega)
    \oplus r^{a-2}H^{m + 3/2}_b(\pa \Omega),
\end{equation*}
for all $|a| < \eta$.
\end{theorem}

See \cite{BNZ1, Kondratiev67, MitreaNistor} or \cite{NP}. If
$\Omega$ is a polygon with maximum angle $\alpha_{M}$, then we can
choose $\eta = \pi/\alpha_M$.

For a convex polytope $\Omega$, \mate{K} will be an integral
operator in a distinguished class of pseudodifferential operators
on the boundary \mate{\pa \Omega}, a class closely related to
\mate{\Diff_{\VV} (\pa\Omega)}, for a suitable Lie algebra of
vector fields \mate{\VV} on \mate{\pa \Omega}. In Melrose's
terminology, this class of pseudodifferential operators
``quantizes'' \mate{\Diff_{\VV} (\pa\Omega)}. These operators can
be thought of as ``singular pseudodifferential operators on
\mate{\pa \Omega}.'' One is lead therefore to consider the
following problem, which we have dubbed ``Melrose's quantization
problem,'' \cite{meicm}:

\begin{quote}
{\bf Melrose's quantization problem:}\ {\em Given a Lie manifold
\mate{(M, \VV)}, one wants to construct \mate{\Psi_{\VV}^\infty(M)},
an algebra of pseudodifferential operators on $M$ with the symbolic
and analytic properties similar to those of the algebra of
pseudodifferential operators on a compact manifolds and such that all
differential operators in \mate{\Psi_{\VV}^\infty(M)} be generated by
$\VV$.}
\end{quote}

If \mate{M = \pa \Omega}, then a variant of the algebra
\mate{\Psi_{\VV}^\infty(M)} should contain the operator \mate{K}
and be compatible with (i.e. quantize) \mate{\Diff_{\VV}(\pa
\Omega)}, thus generalizing the Hardy type operators and the
$b$-calculus. Below, we shall give a construction of the Lie
manifold  $(M, \VV)$ associated to a convex polytope $\Omega$.

\section{Melrose's quantization problem}

We propose a geometric solution in Melrose's spirit. This solution,
given in \cite{aln2}, requires the choice of an appropriate metric on
\mate{M_0 = M \smallsetminus \pa M}, the interior of $M$. More
precisely, we choose on \mate{M_0 = M \smallsetminus \pa M} a metric
\mate{g_0} that has in a neighborhood of any point \mate{x \in M}, a
local orthonormal basis given by sections of \mate{\VV}. A metric
$g_0$ on $M_0$ with this property will be called {\em compatible} (with
the Lie manifold structure $(M, \VV)$). For points \mate{x \in M_0},
the above condition defining a compatible metric is automatically
satisfied, as it follows from Condition (iv) of Definition
\ref{def.LieM}.

The definition of a compatible metric on $M_0$ can be reformulated
as follows. Let $A \to M$ be the Lie algebroid of $(M, \VV)$, that
is, the vector bundle such that $\Gamma(A) \simeq \VV$ as
$\CI(M)$-modules. See our discussion after Definition
\ref{def.LieM}. Then any metric on $A$ defines, by restriction, a
metric on $TM_0$. The resulting metric $g_0$ on $M_0$ is a
compatible metric, and any compatible metric arises in this way.
The metric \mate{g_0} is {\em not} the restriction of a smooth
metric on $M$, in fact, $g_0$ will be singular on \mate{M}.

Let
\begin{equation*}
    (x, y) \mapsto (x, \tau(x, y)) \in TM_0
\end{equation*}
be a local inverse of the Riemannian exponential map \mate{TM_0
\ni v \mapsto \exp_x(-v) \in M_0 \times M_0}. Let
\begin{equation*}
     \big[a_\chi(D)u\big](x) = (2\pi )^{-n}
    {\int_{M_{0}} \left
    (\int_{T^{*}_{x}M_{0}} e^{i \tau(x,y) \cdot \eta} \chi(x,
    \tau(x,y)) a(x, \eta)u(y)\, d \eta \right) dy.}
\end{equation*}
Let $S^m(A^*)$ denote the space of symbols of type $(1, 0)$ (\ie
satisfying H\"ormander's usual estimates \cite{hor3}).

\begin{definition}\label{def.PSM}\
We define \mate{\Psi_{\VV}^{\infty}(M)} to be the space of
pseudodifferential operators $\CIc(M_0) \to \CIc(M_0)$ linearly
generated by \mate{a_\chi(D)} and \mate{b_\chi(D)\exp(X_1) \ldots
\exp(X_k)}, where \mate{a \in S^\infty(A^*)}, \mate{b \in
S^{-\infty}(A^*)}, and \mate{X_j \in \VV}.
\end{definition}

The above definition is consistent with the general principle that all
quantities on \mate{M_0} (functions, Sobolev spaces) should be defined
using the metric \mate{g_0}. See also \cite{MelroseScattering}. This
is, in fact, what leads to the definition of the spaces
\mate{H^m_b(\Omega)} and \mate{H^m_b(\pa \Omega)} before Theorem
\ref{thm.reg1}.

Melrose's quantization problem has a solution \cite{aln2} (see also
\cite{LN1, Monthubert, NWX}. Important related results were obtained
by \cite{emm91, Lauter, LauterMoro, Mazzeo, MelroseScattering,
MelroseMendoza, ScSC, Schulze98, jaredduke}. Let $\Diff(M_0)$ denote
all differential operators on $M_0$.

\begin{theorem}[Ammann-Lauter-Nistor]\label{thm.quant}\
The space \mate{\Psi^{\infty}_\VV(M)} is an algebra of
pseudodifferential operators that ``quantizes'' the Lie algebra
\mate{\VV}, in the sense that \mate{\Psi^{\infty}_\VV(M)} has the
usual symbolic and analytic properties that pseudodifferential
operators have on compact manifolds, and \mate{\Psi^{\infty}_\VV(M) \cap
\Diff(M_0) = \Diff_{\VV}(M)}. In particular, there exist
surjective principal symbol maps $\sigma^{(m)} : \PSM{m} \to
S^m(A^*)$ with kernel $\Psi^{m-1}_\VV(M)$ and any $P \in \PSM{m}$
defines a continuous map $H^s(M_0) \to H^{s-m}(M_0)$.
\end{theorem}

By slightly enlarging the construction of the algebra
$\Psi^{m}_\VV(M)$ by including some additional regularizing
operators, we recover the Hardy type operators as well as the
(small) $b$-calculus.

The most difficult part in the proof of the above theorem is to show
that $\PSM{\infty}$ is closed under composition. Our proof in
\cite{aln2} is to show that $\PSM{\infty}$ is the homomorphic image of
$\Psi^\infty(\maG)$. Here $\maG$ is a groupoid integrating the Lie
algebroid $A$ associated to $M$ \cite{CrainicFernandez, NistorINT} and
$\Psi^\infty(\maG)$ is algebra of pseudodifferential operators
\cite{NWX} (in particular, it is closed under composition). See also
\cite{connes79, Monthubert}. The groupoid $\maG$ plays the role of a
kernel space, because $\Psi^\infty(\maG) = I_c^\infty(\maG, M)$, the
space of compactly supported distributions on $\maG$ that are conormal
to $M$. These kernel spaces are very closely related to a construction
of Melrose (the stretched $b$-product ${}^bM^2$), see \cite{meicm,
MelroseScattering}. It is, in general, a difficult task to find a
groupoid integrating a Lie algebroid $A$, and, in fact, this is not
always possible. It is a deep theorem of Crainic and Fernandez that
this is possible for the Lie algebroids associated to Lie manifolds
\cite{CrainicFernandez}. See also \cite{NistorINT}, which suffices for
example for the examples considered in next section.

\section{An application to Fredholm conditions\label{sec.FC}}

We shall now obtain some criteria for operators $P \in \PSM{m}$ to be
Fredholm. This has applications to boundary value problems as well as
to non-linear partial differential equations on non-compact manifolds.

We define the Sobolev space $H^s(M_0)$ to be the domain of $(1 +
d^*d)^{s/2}$, where $d$ is the de Rham differential and $s \ge 0$. For
$s < 0$, we use duality to define $H^s(M_0)$. Also, an elliptic
operator $P \in \PSM{m}$ is one for which $\sigma^{(m)}(P)(\xi) \neq
0$ for all $\xi \in A^*$, $\xi \neq 0$.

The main theorem, which will ocupy us the rest of this section, is the
following (``type I Lie manifolds,'' as well as the rest of the
unexplained notation of the following theorem, are introduced below).

\begin{theorem}\label{thm.Fredholm}\
Let \mate{(M, \VV)} be a type I Lie manifold. Assume each hyperface of
$M$ has a defining function. If \mate{P \in \PSM{m}}, then there exist
pseudodifferential operators \mate{P_\alpha} on \mate{M_\alpha \times
G_\alpha}, invariant with respect to right translations by $G_\alpha$
such that \mate{P : H^{s}(M_0) \to H^{s-m}(M_0)} is Fredholm if, and
only if, \mate{P} is elliptic and all \mate{P_\alpha : H^{s}(M_\alpha
\times G_\alpha) \to H^{s-m}(M_\alpha \times G_\alpha)} are invertible
for all $\alpha \neq 0$.
\end{theorem}

This theorem will follow from the results of \cite{LMN1, LN1}. More
general Fredholmness conditions were obtained in \cite{LMN1, LN1}, but
they involve some conditions that may be difficult to use. On the
other hand, the conditions for the above theorem are easier to check.

We shall assume that our Lie manifold $(M, \VV)$ satisfies the
following four conditions. Our first condition is that there exist,
for any (closed) face $F \subset M$, a fibration $p_F : F \to B_F$
with connected fibers, such that
\begin{equation}\label{eq.C1}
    \varrho(A_p) = T_pp_F^{-1}\big(p_F(p)\big),
    \text{ for any } p \in F_0:=\overset{\circ}{F}.
\end{equation}

We shall use the Lie algebroid $A \to M$ associated to $\VV$ and the
anchor map $\varrho : A \to TM$ introduced in Remark
\ref{rem.Algebroid}.

Another way of formulating our first condition, Equation
\eqref{eq.C1}, is that, for any $p$ in the interior of $F$, the
tangent space at $p$ through the fiber of $p_F$ containing $p$
coincides with the set $X(p)$, $X \in \VV$. Yet another way of
formulating this condition is that the set $\{\exp_X(p)\}$ is the
fiber of $p_F$ containing $p$, where $\exp_{tX}$ is the one-parameter
group of diffeomorphisms obtained by integrating $X$ and $p$ is in
$F_0$, the interior of $F$. (See \cite{aln1} for the easy proof that
$\exp_{tX}$ is defined for all $t$.) From this condition it also
follows that the isotropy Lie algebras
\begin{equation}\label{eq.mfl}
    \mathfrak l_p := \ker (\varrho : A_p \to T_pM),
    \quad p \in F_0,
\end{equation}
have the same dimension, and hence they define a vector bundle on
the interior of $F$. Let $\mfL_F \to F_0$ denote this vector
bundle. Then $\mfL_F$ is a bundle of Lie algebras.

For any $p \in M$, the vector space $\mathfrak l_p$ (see Equation
\eqref{eq.mfl}) has a natural structure of Lie algebra. (To see this,
let $X, Y \in \VV$ be such that $X(p) = Y(p) = 0 \in T_pM$.  Then $[X,
Y](p)$ depends only on $X(p)$ and $Y(p)$. See also \cite{Mackenzie}.)
We shall call $\mathfrak l_p$ {\em the isotropy Lie algebra of
$p$}. Let $G_p$ be the simply connected Lie group with Lie algebra
$\mathfrak l_p$. We say that $G_p$ is an exponential Lie group if the
exponential map defines a diffeomorphism
\begin{equation}\label{eq.pre.C2}
    \exp : \mathfrak l_p \simeq G_p.
\end{equation}
The simply-connected nilpotent Lie groups and most simply-connected
solvable Lie groups satisfy this condition. Our third condition is
then
\begin{equation}\label{eq.C2}
    G_p \text{\em\ is a solvable, exponential group}.
\end{equation}

To formulate the third condition, let $M_\alpha$ be the orbits of
the diffeormorphisms $\exp_{X}$, $X \in \VV$, acting on $M$, where
$\alpha$ belongs to an index set $\mathfrak I$ containing $0$. By
our first assumption, if $p \in M_\alpha$ is an interior point of
a face $F$, then $M_\alpha$ coincides with the interior of the
fiber of $p_F : F \to B_F$ containing $p$. In particular,
$M_\alpha = M_0$ if $\alpha = 0$. Our third assumption is that
there exists a bundle of Lie algebras $\mfA_F$ on the interior of
$B_F$ such that
\begin{equation}\label{eq.C3}
    \mfL_F \simeq p_F^*\mfA_F.
\end{equation}
Let $\mfA_{Fq}$ be the fiber of $\mfA_F$ above $q$, with $q$ in
the interior of $B_F$. Let $G_q$ be a simply-connected Lie group
with Lie algebra $\mfA_{Fq}$. Let $q = p_F(p)$. Since $G_q \simeq
G_p$ and $\mfA_{Fq} \simeq \mathfrak l_p$, it follows from our
second assumption (Equation \eqref{eq.C2}) that the exponential
map
\begin{equation}\label{eq.exp}
    \exp : \mfL_{Fq} \to G_q
\end{equation}
is a diffeomorphism.

We shall also need differentiable groupoids. Let say first that a
groupoid is a ``group with several units,'' and that the product of
two elements is defined only if the domain of the first matches the
range of the second one. The model for a groupoid is a set of
bijective functions. More precisely, a groupoid is a small category
all of whose elements are invertible. (A category is small if the
class of its objects is in fact a set.)  See \cite{LN1, NWX} for an
introduction to differentiable groupoids that is suitable for our
purposes.

Let $p_F : G_F \to \stackrel{\circ}{B_F}$ be $\mfA_F$ as a manifold,
but with the Lie group structure on each fiber induced by the
exponential map (which is a diffeomorphism, see Equation
\eqref{eq.exp}). For each $\alpha$, let $q = p_F(M_\alpha)$. We shall
denote by $G_\alpha$ the fiber of $G_F$ above $q$. (This makes sense
in view of our first assumption, Equation \eqref{eq.C1}.) Consider the
fibered product
\begin{multline}\label{def.G.F}
    \maG_F = F_0 \times_{B_F} F_0 \times_{B_F} G_F \\ := \{(x, y, g)
    \in F_0 \times F_0 \times G_F, p_F(x) = p_F(y) = p_F(g) \in B_F
    \},
\end{multline}
with the groupoid structure given by the product $(x, y, g)(y, z, h) =
(x, z, gh)$, the set of units $F_0$, and the domain map $d(x, y, g) =
y \in F_0$ and range map $r(x, y, g) = y \in F_0$. Then $\maG_F$ is a
differentiable groupoid.  Our fourth assumption is that
\begin{equation}\label{eq.C4}
    \maG := \cup \maG_F \,\text{\em is a Hausdorf differentiable
    groupoid with Lie algebroid } \simeq A.
\end{equation}
The groupoid $\maG$ is the disjoint union of the groupoids $\maG_F$,
and the structural morphisms (composition, domain, range, ...) are the
ones induced from $\maG_F$. In particular, the set of units of $\maG$
is the disjoint union of the units of the groupoids $\maG_F$, that is,
$\maG$ has as a set of units $M = \cup F_0$. By the results of
\cite{NistorINT}, there is at most one differentiable structure on
$\maG$ with Lie algebroid $A$. It induces the given differentiable
structure on $\maG_F$.

\begin{definition}\label{def.typeI}\
A Lie manifold satisfying the above four conditions (Equations
\eqref{eq.C1}, \eqref{eq.C2}, \eqref{eq.C3}, and \eqref{eq.C4})
will be called {\em a type I Lie manifold}.
\end{definition}

For the proof of Theorem \ref{thm.Fredholm} it is necessary to
recall a few constructions and to prove some intermediate results.

Recall \cite{LN1, NWX} that $P \in \Psi^{\infty}(\maG)$ is in fact a
family of pseudodifferential operators $P = (P_x)$, $x \in M$, with
$P_x$ acting on
\begin{equation*}
    \maG_x := d^{-1}(x) = \{(z, x, g), z \in M_\alpha,\ e \text{
    the unit of } G_\alpha \} \simeq M_\alpha \times G_\alpha
\end{equation*}
and satisfying some additional assumptions. These additional
conditions are:\ right invariance for multiplication by elements in
$\maG$, that the family $(P_x)$ be a smooth family, and a suport
condition that implies, in particular, that all $P_x$ are properly
supported.

If $x, y \in M_\alpha$, then right multiplication by $(x, y, e)$ is a
diffeomorphism $\maG_x \to \maG_y$ mapping $P_x$ to $P_y$, by the
assumption of right invariance. The restriction $\tilde P_\alpha$ of
$P$ to $\maG_x$ is simply $P_x$, where $x \in M_\alpha$.  The
canonical isomorphism $\maG_x \to M_\alpha \times G_\alpha$ will map
\begin{equation}\label{eq.defPa}\
    P_x \to P_\alpha
\end{equation}
for all operators $P_x$, with $P_\alpha$ a pseudodifferential operator
on $M_\alpha \times G_\alpha$ that is independent of $x$. This
operator is right invariant with respect to the action of $G_\alpha$,
again by the invariance condition.

Let $\phi \in \CIc(\maG)$.  We shall denote by $\phi_x$ the
restriction of $\phi$ to $\maG_x$. We fix a metric on $A$ which will
fix a metric on each of the spaces $\maG_x$ and hence a volume form
smoothly depending on $x$. We shall denote by $\| \,\cdot\, \|$ the
norm on $L^2(\maG_x)$ or the norm of a bounded opearator on this
space, for any $x$. There will be no danger of confusion. We begin by
examining the consequences of the assumption that $\maG$ is Hausdorf.

\begin{proposition}\label{prop.lsc}\
Let $P = (P_x) \in \Psi^{m}(\maG)$. Then $\|P_x \phi_x\|$ depends
continuously on $x \in M$, for any $\phi \in \CIc(\maG)$.
\end{proposition}

\begin{proof}\
Fix $\phi \in \CIc(\maG)$.  From the definition of the algebra
$\Psi^\infty(\maG)$, it follows that there exists a function $\psi \in
\CIc(\maG)$ such that $\psi_x = P_x \phi_x$. The continuity of the
function $\|\psi_x\|$ follows from the assumption that $\maG$ is
Hausdorf and from the smoooth dependence of the measure on $\maG_x$ on
$x$.
\end{proof}

We now prove as a consequence the following result.

\begin{corollary}\label{cor.inj}\
Let $P = (P_x) \in \Psi^{m}(\maG)$ and $x \in M_0$. If $P_x = 0$, then
$P_y = 0$, for any $y \in M$. That is $P = 0$.
\end{corollary}

\begin{proof}\
This is a consequence of the fourth Assumption, namely Equation
\eqref{eq.C4}. Indeed, assume $P_x = 0$ for some $x \in M_0$. Then
$P_y = 0$ for all $y \in M_0$, by the right invariance of the
operators $P_x$. To prove that $P_y = 0$ for some arbitrary $y$, we
now show that $P_y \eta = 0$ for any $\eta \in \CIc(\maG_y)$.  Let
$\phi \in \CIc(\maG)$ that restricts to $\eta$ on $\maG_y$ (\ie $\phi_y
= \eta$). This is possible since $\maG_y$ is a closed subset of the
Hausdorf, locally compact space $\maG$. Then $\|P_y \phi_y\|$ is a
continuous function of $y \in M$ that vanishes for $y \in M_0$. Since
$M_0$ is dense in $M$, we obtain that $P_y \phi_y = 0$ for all $y$.
\end{proof}

The assumption that $\maG$ is Hausdorff therefore implies that the
natural action of $\Psi^{\infty}(\maG)$ on $\CIc(M_0)$ is faithful
(\ie the induced morphism $\Psi^{\infty}(\maG) \to
\operatorname{End}(\CIc(M_0))$ is injective). Fix $z \in M_0$ and
consider the canonical bijection (diffeomorphism) $\maG_z \to
M_0$. Then the map
\begin{equation}
    \Psi^{\infty}(\maG) \ni P = (P_x) \to P_z \in
    \Psi^{\infty}(\maG_z) \simeq \PSM{\infty}
\end{equation}
is a bijection. We shall henceforth identify these two algebras (this
is incidentaly the canonical surjection constructed in
\cite{aln2}). In particular, we can define $P_\alpha$, for any $P \in
\PSM{\infty} = \Psi^{\infty}(\maG)$ using Equation \eqref{eq.defPa}.

\begin{corollary}\label{cor.mult}\
We have that $(PQ)_\alpha = P_\alpha Q_\alpha$, for all $P, Q \in
\PSM{\infty}$.
\end{corollary}

\begin{proof}\
The product in the algebra $\Psi^\infty(\maG)$ is $PQ = (P_x Q_x)$, if
$P = (P_x)$ and $Q = (Q_x)$, $x \in M$. The result then follows from
the definition of $P_\alpha$ given in Equation \eqref{eq.defPa}.
\end{proof}

Yet another corollary of Proposition \ref{prop.lsc} is the following.

\begin{corollary}\label{cor.lsc}\
Let $P = (P_x) \in \Psi^{0}(\maG)$. Then the
function $M \ni x \to \|P_x\| \in \RR$ is lower
semi-continuous.
\end{corollary}

\begin{proof}
Indeed, let $\alpha \in \RR$. We need to show that the set $\{x \in M,
\|P_x\| > \alpha \}$ is open in $M$. Let $y \in M$ be such that
$\|P_y\| > \alpha$. Then we can find $\eta \in \CIc(\maG_y)$ such that
$\|P_y \eta\| > \alpha \|\eta\|$.  Let $\phi \in \CIc(\maG)$ be such
that $\phi_y = \eta$. Since $\|P_x \phi_x\|$ and $\|\phi_x\|$ are
continuous, we have that $\|P_x \phi_x\|/\|\phi_x\|$ is well defined
and continuous in a neighborhood of $y$. But then $\|P_x
\phi_x\|/\|\phi_x\|>\alpha$ defines an open neighborhood of $y$ on
which $\|P_x\| > \alpha$.
\end{proof}

This in turn gives the following.

\begin{corollary}\label{cor.B.inj}\
For any $P_x \in \Psi^0(\maG)$ we have $\|P_y\| \le \|P_x\|$ for any
$x \in M_0$, $y \in M$. In other words, the function $M \ni y \to
\|P_y\| \in [0, \infty)$ attains its maximum at any point $x \in M_0$.
\end{corollary}

\begin{proof}\
All operators $P_x$ are unitarily equivalent for $x \in
M_0$. Therefore the function $M \ni y \to \|P_y\| \in [0, \infty)$ is
constant on $M_0$. Now if $\|P_y\| > \|P_x\|$ for some $y \in \pa M =
M \smallsetminus M_0$ and some $x \in M_0$, then, by choosing $\|P_y\|
> \alpha > \|P_x\|$, we contradict the fact that the set $\{y \in M,
\|P_y\| > \alpha \}$ is open in $M$.
\end{proof}

Let $\overline{\Psi}_{-\infty}$ be the closure of the ideal
$\PSM{-\infty}$ in the family of norms of operators $H^{t}(\maG_x)
\to H^{r}(\maG_x)$, $x \in M$. By Corollary \ref{cor.B.inj}, this
closure is the same as the closure of $\overline{\Psi}_{-\infty}$
in the topology of continuous operators $H^{t}(M_0) \to
H^{r}(M_0)$.  Let $\Psi^{s} := \PSM{s} + \Psi^{-\infty}$.  Then
$\overline{\Psi}_{s}\overline{\Psi}_{s'} \subset
\overline{\Psi}_{s+s'}$.

Denote by $\maL(\maH)$ the set of continuous linear operators
$\maH \to \maH$. Let us notice that the exact sequence of
envelopping $C^*$-algebras of groupoids (see for example
\cite{LN1}[Equation 16]) or the structure theorem
\cite{LN1}[Theorem 4.4] show that the natural representation
$C^*(\maG) \to \maL(\maH)$ is injective. In particular, $\maG$ is
amenable. Then Theorems \ref{thm.Fredholm} and \ref{thm.comp}
follow right away from \cite{LN1}[Theorem 9.]. We prefer however
to include some arguments, to make the paper more complete. Let
$\mathfrak A(M)$ be the norm closure of $\PSM{0}$ acting on
$L^2(M_0)$. The Theorems \ref{thm.Fredholm} and \ref{thm.comp}
remain true for $P \in \mathfrak A(M)$ and $m = 0$.

We shall need also the following theorem.

\begin{theorem}\label{thm.comp}\
We keep the assumptions and notation of Theorem \ref{thm.Fredholm}
and fix $s \in \RR$. Let \mate{P \in \PSM{m}}, then \mate{P :
H^{s}(M_0) \to H^{s-m}(M_0)} is compact if, and only if,
\mate{\sigma^{(m)}(P) = 0} and \mate{P_\alpha = 0}, for all
$\alpha \neq 0$.
\end{theorem}

\begin{proof}\
This will be a consequence of the results of \cite{LMN1, LN1}.  As in
\cite{alnv}[Proposition 5.2 and Theorem 6.2], we can find an
invertible pseudodifferential operator $P \in \overline{\Psi}_{r}$,
for any $r$. (``Invertible'' here means that the inverse is in
$\overline{\Psi}_{-r}$.) This allows us to assume that $P$ has order
zero.

Let $x_H$ be a defining function for each hyperface $H$ of $M$ and $x$
the product of all defining functions of hyperfaces of $M$.  Let $P
\in \PSM{-1}$. If $P_\alpha = 0$ for all $\alpha \neq 0$, then $P =
xQ$, with $Q \in \PSM{-1}$. Therefore $P$ maps $H^s(M_0)$ continuously
to $xH^{s-1}(M_0)$. Since $xH^{s-1}(M_0) \to H^s(M_0)$ is a compact
map (see, for example, \cite{AIN}[Theorem 3.6] for this easy
generalization of Kondrachov's theorem), it follows that $P : H^s(M_0)
\to H^s(M_0)$ is compact.

As above, we can assume that $P$ has order zero. Suppose now that
$P : H^s(M_0) \to H^s(M_0)$ is compact. Then $\sigma^{(0)}(P) =
0$, as in the classical case \cite{hor3}. Assume, by
contradiction, that $P_\alpha \neq 0$, for some $\alpha$. Fix for
the rest of this discussion $x \in M_\alpha$ and $\phi \in
\CIc(\maG_x)$, $\maG_x = M_\alpha \times G_\alpha$, such that
$P\phi \neq 0$. We extend $\phi$ to a smooth, compactly supported
function on $\maG$, still denoted by $\phi$. Let $\phi_y$, $y \in
M_0$ be the restriction of $\phi$ to $\maG_y \simeq M_0$. As $y
\to x$, $y \in M_0$, we have that $\phi_y \to 0$ weakly, but
$\|P\phi_y\| \to \|P_x \phi_x\| \neq 0$. So $P$ cannot be compact.
\end{proof}

We shall need also the following corollary of the above proof.

\begin{corollary}\label{cor.proof}\ We keep the notation of the
proof of Theorem \ref{thm.comp}.

\noindent {\rm (i)}\ $(PQ)_\alpha = P_\alpha Q_\alpha$ for $P, Q
\in \overline{\Psi}_{\infty}$.

\noindent {\rm (ii)}\ Assume that $P \in \PSM{m}$ is a Fredholm
operator \mate{P : H^{s}(M_0) \to H^{s-m}(M_0)}. Then there exist
$Q \in \overline{\Psi}_{-m}$ such that $P Q - I$ and $Q P - I$ are
compact operators.
\end{corollary}

\begin{proof}\
Part (i) is clear by the definition of $\overline{\Psi}_{\infty}$ and
Corollary \ref{cor.B.inj}.

As in the above proof, we can assume that $P$ has order zero. It was
proved in \cite{alnv} and in \cite{LMN1} that $\overline{\Psi}_{0}$ is
closed under holomorphic functional calculus. Since we can construct
$Q$ out of $P$ using holomorphic functional calculus, it follows that
$Q \in \overline{\Psi}_{0}$.
\end{proof}

We are ready now to prove the main result of this section, Theorem
\ref{thm.Fredholm}.

\begin{proof}\
Assume that $\sigma^{(0)}(P)(\xi)$, $\xi \neq 0$, and
$P_\alpha$, $\alpha \neq 0$, are invertible.  The structure theorems
of \cite{LMN1, LN1} show that the map $\overline{\sigma}$
\begin{equation}
    \PSM{0} \ni P \to \big(\sigma^{(0)}(P)\vert_{S^*A},
    P_\alpha) \in C(S^*A) \oplus \oplus_\alpha
    \maL(L^2(M_\alpha \times G_\alpha))
\end{equation}
extends to $\mathfrak A(M)$, the the norm closure of $\PSM{0}$
acting on $L^2(M_0)$. Moreover, the structure theorems of
\cite{LN1}[Theorem 4.4] (see also \cite{LMN1}) also show that the
kernel of the map $\overline{\sigma}$ is given by the set of
compact operators. (We are using here also the fact that solvable
groups are amenable and hence that any irreducible
$*$-representation of $\PSM{0}$ is contained in one of the
representations on $L^2(M_\alpha \times G_\alpha)$.) Therefore $P$
is invertible modulo compact operators if, and only if,
$\overline{\sigma}(P)$ is invertible. This completes the proof of
Theorem \ref{thm.Fredholm}.
\end{proof}

The first part of the above theorem has an elementary proof as
follows. Choose $Q \in \overline{\Psi}_{0}$ such that $PQ - I$ and
$QP - I$ are compact, using Corollary \ref{cor.proof}. Then
$P_\alpha Q_\alpha - I = (PQ - I)_\alpha = 0$, by Theorem
\ref{thm.comp}. Similary, $Q_\alpha P_\alpha - I = 0$. This proves
that $P_\alpha$ is invertible, for all $\alpha \neq 0$. The
ellipticity of $P$ follows from classical results \cite{hor3}.
(See \cite{MelroseNistor1} for the details of this argument.) An
elementary proof of the second part of the above theorem is
usually obtained by constructing geometrically a bounded right $s$
inverse of $\overline{\sigma}$ ($s$ is defined on the range of
$\overline{\sigma}$).

Earlier related results were obtained by \cite{CO2, Kondratiev67,
Lauter, LauterMoroianu, Mazya, Mazzeo, MelroseScattering, Mendoza,
SchroheSI, SchroheFC, Schulze98}.

It is interesting to notice that each of the operators
\mate{P_\alpha} is \mate{G_{\alpha}}--invariant and ``of the same
kind'' as the operator \mate{P}, for example, if $P$ is the
Laplace operator associated to a compatible metric $g$, then
$P_\alpha$ will be the Laplace operator corresponding to the
induced metric on $M_\alpha \times G_\alpha$. See \cite{aln1} for
more results in this direction. This leads to an inductive
procedure to study an operator \mate{P \in \PSM{\infty}}, which
will be used, for example, in Section \ref{sec.Sp}.

\section{Examples\label{sec.Ex}}

Let us discuss some examples of how the above theory can be used
to study concrete examples. Most of these examples go back to
Melrose \cite{MelroseScattering}. Since $M_0$ is always the
interior of $M$ and the group $G_0$ is reduced to only one
element, we shall typically assume below that $\alpha \neq 0$.

\begin{example}\label{ex.one}\
Let \mate{M} be a manifold with smooth, connected boundary
\mate{\pa M}, \mate{M_0 = M \smallsetminus \pa M}, as before. On
$M$ we consider the set \mate{\VV = \VV_b} of vector fields that
are {\em tangent} to \mate{\pa M}. We impose no condition on these
vector fields in the interior, as required by Axiom (iv) of the
definition of a Lie manifold, Definition \ref{def.LieM}. Let
\mate{y_2, \ldots, y_n} be some local coordinates on \mate{\pa M}
and let $x$ denote the distance to the boundary. At the boundary
\mate{\pa M = \{x = 0\}}, a local basis of $\VV_b$ is given by
${x}\pa_x, \pa_{y_2}, \ldots, \pa_{y_n}$.

An example of a compatible metric on $M_0$ is \mate{g_0 =
\frac{(dx)^2}{x^2} + h}, with \mate{h} smooth on \mate{M}. The
resulting algebra \mate{\Diff_{\VV_b}(M)} of differential
operators is the algebra of totally characteristic differential
operators. The metric on \mate{M_0} is that of a manifold with
cylindrical ends. The resulting pseudodifferential calculus is the
subalgebra of properly supported pseudodifferential operators in
Melrose's $b$-calculus $\Psi_{b}(M)$. The Lie algebroid $A \to M$
is Melrose's compressed tangent bundle ${}^bTM$. The groupoid
integrating ${}^bTM$ is obtained from Melrose's stretched
$b$-product by removing the faces not intersecting the diagonal.

In this example \mate{\{\alpha \neq 0\} = \{ \pa M \}} consists of
exactly one element and  \mate{M_\alpha = \pa M}, \mate{G_\alpha =
\RR}. The Fredholmness criteria were obtained in increasing
generality in \cite{Kondratiev67, LockhartOwen, MelroseMendoza}.
\end{example}

This example is basic in that it helps us understand easier other,
more complicated examples. In the following examples we will
indicate only what is different from the first example.

\begin{example}\
Take now \mate{\VV_0} to be the space of vector fields on \mate{M}
that vanish on \mate{\pa M}. At the boundary \mate{\pa M = \{x =
0\}} a local basis is given by \mate{ \rx \pa_x,  \rx \pa_{y_2},
\dots,  \rx \pa_{y_n}}. The resulting geometry is that of an
asymptotically hyperbolic manifold. The strata different from
$M_0$ are \mate{M_\alpha = \{\alpha\}} are parametrized by
\mate{\alpha \in \pa M}. The group $G_{\alpha} = T_{\alpha}(\pa M)
\rtimes \RR$ is a solvable Lie group with $t \in \RR$ acting by
dilation by $e^t$ on $T_{\alpha}(\pa M)$.

Recently these manifolds have been used in Mathematical physics in
connection to the AdS--CFT correspondence \cite{anderson:p03, Lee,
MazzeoQG, GraZwor}. Earlier, slightly larger larger algebras of
pseudodifferential operators quantizing $\VV_0$ were constructed
in \cite{Mazzeo, Schulze98} and called the ``edge-calculus.''
\end{example}

We now discuss an example that generalizes the manifolds Euclidean
at infinity.

\begin{example}\
Let us take now \mate{\VV_{sc}} to be the space of vector fields
on \mate{M} that vanish on \mate{\pa M} and have the property that
their normal component to the boundary vanishes of second order at
the boundary. At the boundary \mate{\pa M = \{x = 0\}} a local
basis is given by \mate{ {x^2} \pa_x, \rx \pa_{y_2}, \ldots,  \rx
\pa_{y_n}}. The resulting geometry is that of an asymptotically
flat manifold. As in the previous example, \mate{\{\alpha \neq
0\}= \pa M}, each \mate{M_\alpha = G_{\alpha} = T_{\alpha}(\pa M)
\times \RR} is an {abelian} Lie group, and each \mate{P_\alpha} is
\mate{G_\alpha} invariant.

This example is the best understood so far. For example, earlier
versions of the pseudodifferential calculus were introduced by
Parenti (called the ``SG-calculus'') \cite{Parenti} and Melrose
\cite{MelroseScattering} (called the ``scattering-calculus''). See
\cite{Ashtekar} for an application of asymptotically flat
manifolds to Quantum Gravity.
\end{example}

Here is now an example similar to that of asymptotically
hyperbolic manifolds considered above. This example is relevant
for the analysis on locally symmetric spaces and for boundary
value problems on polyhedral domains.

\begin{example}\
Let \mate{\pi: \pa M \to B} be a fibration, and let
\mate{\VV_{\pi}} be the space of vector fields on \mate{M} that
are tangent to the fibers of this fibration. We choose a system of
coordinates at the boundary \mate{\pa M = \{x = 0\}} such that the
fibration becomes a product in that neighborhood. Then a local
basis of $\VV_{\pi}$ on the domain of our coordinate chart is
given by \mate{\rx \pa_x, \rx \pa_{y_2}, \ldots , \rx \pa_{y_k},
\pa_{y_{k + 1}}, \ldots, \pa_{y_{n}}}.

In this example, the set of non-zero parameters is \mate{\{\alpha
\neq 0\}= B}, the strata is given by \mate{M_\alpha =
\pi^{-1}(\alpha), \alpha \in B}, and \mate{G_{\alpha} =
T_{\alpha}B \rtimes \RR} is a solvable Lie group with $\RR$ acting
again by dilations. Earlier, slightly larger larger algebras of
pseudodifferential operators quantizing $\VV_0$ were constructed
in \cite{Mazzeo, Schulze98} and called the ``edge-calculus.''
\end{example}

We now include an example of a Lie manifold that is not type I. It
a variation of the previous example. It is not clear how to
generalize Theorem \ref{thm.Fredholm} to this example, although
Fredholmness conditions can be obtained as in \cite{LMN1}.

\begin{example}\
Let \mate{F \subset T \pa M} be a foliation of the boundary of $M$. We
assume that not all leaves of $F$ are closed in $M$, to avoid
trivialities.  We take then \mate{\VV = \VV_F} to be the space of
vector fields on \mate{M} that are tangent to the leaves of
\mate{F}. No earlier pseudodifferential calculi on these manifolds
were considered before.
\end{example}

We conclude with an example that generalizes our first example to
manifolds with corners.

\begin{example}\label{ex.cor}\
Let \mate{M} be a compact manifold with corners. We define
\mate{\VV = \VV_b} to be the space of vector fields on $M$ that
are tangent to all hyperfaces of \mate{M}, \cite{MelroseNistor1,
MelrosePiazza}. In this example, \mate{\{\alpha \neq 0\}} is the
set of faces \mate{H} of maximal dimension of \mate{M} (\ie the
hyperfaces of $M$) and \mate{M_H = H} for any hyperface $H$.
Finally, \mate{G_\alpha = \RR}. See also \cite{MelroseNistor1}. A
Riemannian manifold isometric to $M_0$ with a compatible metric is
called {\em a manifold with multi-cylindrical ends}.
\end{example}

We now discuss the Lie manifold with boundary associated to a
convex polytope. They are type I.

\begin{example}\ Let $\PP$ be a simplex in $\RR^N$. Let
$(\Sigma(\PP), \kappa)$ be its desingularization, where
$(\Sigma(\PP), \VV)$ is a Lie manifold with boundary, as in
\cite{AIN}. Then $\pa \Sigma(\PP)$ and the double of $\Sigma(\PP)$
are type I Lie manifolds. (The ``double'' of $\Sigma(\PP)$ is
obtained by gluing two copies of $\Sigma(\PP)$ along their true
boundary, \ie along the closure of the set of boundary points that
correspond to each other and are not at infinity.)
\end{example}

\section{Spectra\label{sec.Sp}}

In this section we give an application of the Fredholmness conditions
to the determination of the spectrum of the Dirac and Laplace
operators on the manifolds arising in Example \ref{ex.cor}. In this
section, we shall assume that $M$ is a manifold with $\pa M \neq
\emptyset$.

Let us consider for a moment the framework of \ref{ex.one}, which
is a particular case of Example \ref{ex.cor}. Let \mate{P =
\Delta_{M_0} - \lambda}. Then
\begin{equation*}
    {P_\alpha = \Delta_{\pa M \times \RR} =
    \Delta_{\pa M} - \pa_t^2 - \lambda}.
\end{equation*}
Let \mate{\hat{P}(\tau) = \Delta_{\pa M} + \tau^2 - \lambda}, be
the Fourier transform of \mate{P_\alpha} in the $t$ variable. This
is what Melrose calls the ``indicial family'' associated to $P$.
Since the spectrum of \mate{\Delta_{\pa M}} is
\begin{equation*}
    { \sigma(\Delta_{\pa M}) = \{ 0, \lambda_1,
    \lambda_2, \ldots \} \subset [0, \infty) },
\end{equation*}
we obtain that \mate{\hat{P}(\tau)} is invertible for any
\mate{\tau \in \RR} if, and only if, \mate{\lambda < 0}. Hence
{$\Delta_{M_0} - \lambda$} is Fredholm, if, and only if,
\mate{\lambda < 0}. This shows that {$\sigma_e(\Delta_{M_0}) = [0,
\infty)$}. But then
\begin{equation*}
    {[0, \infty) \subset \sigma_e(\Delta_{M_0})
    \subset \sigma(\Delta_{M_0}) \subset [0,\infty)}
\end{equation*}
and hence %\begin{equation*}
$\sigma(\Delta_{M_0}) = [0,\infty)$.
%\end{equation*}

This argument generalizes to higher rank spaces \cite{LN1} to
prove the following result that was formulated as a conjectured in
\cite{MelroseScattering}.

\begin{theorem}[Lauter-Nistor]\label{thm.LN}\
Assume \mate{M} is as in Example \ref{ex.cor} and $\pa M \neq
\emptyset$. Then
\begin{equation*}
    {\sigma(\Delta_{M_0}) = [0,\infty)}.
\end{equation*}
\end{theorem}

We now extend the reasoning of the proof of the above theorem in
\cite{LN1} to study the Dirac operator. Recall that in this section
we assume that $\pa M \neq \emptyset$.

\begin{theorem}\label{Theorem.c}\
Let \mate{M} is as in Theorem \ref{thm.LN} and $W \to M$ be a
Clifford bundle over $A^*$.  We assume no face of $M$ has
dimension zero. Let $D_F$ be the Dirac operator on $F$ with
coefficients in $W\vert_F$ for any face $F \subset M$. We assume
that $\ker(D_F) = 0$, for any $F \neq M$. Then each $D_F$ is
invertible and
\begin{equation*}
    \sigma_e(D_M) = (-\infty, -c] \cup [c, \infty).
\end{equation*}
where $c^{-1} = \max \{ \|D_H^{-1}\| \} > 0$, for $H$ ranging
through the set of hyperfaces of $M$. Moreover, $D_M$ is
invertible if, and only if, $\ker(D_M) = 0$.
\end{theorem}

\begin{proof}\
We shall prove this by induction. If $M$ has no boundary, then
$\sigma_e(D_M) = \emptyset$ and $D_M$ is invertible if, and only
if, there exist no $L^2$-harmonic spinors (\ie $\ker(D_M) \neq
0$).  This situation is excluded by our theorem since $\pa M \neq
\emptyset$; it is needed, however, for the inductive hypothesis.
If $M$ has boundary, then $\sigma_e(D_M)$ is the spectrum of $D_1
:= D_{\pa M} + c(dt)\pa_t$ acting on $L^2(\pa M \times \RR, W_F)$,
where $t$ denotes the $\RR$-component and $c(\omega)$ is the
operator of Clifford multiplication by $\omega$. We have $D_1^2 =
-\pa_t^2 + D_{\pa M}^2$. Therefore $\sigma(D_1^2) = [c^2,
\infty)$, where $c^{-1} = \|D_{\pa M}^{-1}\|$, is defined since
there are no $L^2$-harmonic spinors on $\pa M$.

Let $V : L^2(\pa M \times \RR, W_F) \to  L^2(\pa M \times \RR,
W_F)$ be given by $V(u)(t) = c(dt)u(-t)$. Then $V D_1 V^{-1} = -
D_1$, and hence $\sigma(D_1)$ is symmetric with respect to $0$.
Hence $\sigma(D_1) = (-\infty, -c] \cup [c, \infty)$. This proves
that $\sigma_e(D_M) = (-\infty, -c] \cup [c, \infty)$. In
particular, since $0 \not \in \sigma_e(D_M)$, we have that $0 \in
\sigma(D_M)$ if, and only if, $0$ is an eigenvalue of $D_M$.

The inductive step, in general, follows as exactly as in the case
of a manifold with boundary, but replacing $\pa M$ with a
hyperface $H$ of $M$.
\end{proof}

A similar reasoning gives the following.

\begin{theorem}\
We keep the same assumptions as in Theorem \ref{Theorem.c}, except
that we assume that $\ker(D_F) = 0$ for at least one face $F \neq
M$. Then
\begin{equation*}
    \sigma_e(D_M) = \RR.
\end{equation*}
\end{theorem}

\begin{proof}\
Let $D_1$ be the restriction of $D_M$ to the groupoid
corresponding to the face $F$. Then $\sigma(D_1) = \RR$, as in the
proof of Theorem \ref{Theorem.c}. But $\sigma(D_1) \subset
\sigma_e(D_M)$, by Theorem \ref{thm.Fredholm} (or by \cite{LN1,
MelrosePiazza}).
\end{proof}

We obtain the following corollary.

\begin{corollary}\label{cor.c}\
We continue to assume that $M$ has no faces of dimension zero and
keep the same notation as in Theorem \ref{Theorem.c}. Then $D_M$
is Fredholm if, and only if, $D_F$ has no $L^2$-harmonic spinors,
for any face $F \subset M$, $F \neq M$. Similarly, $D_M$ is
invertible if, and only if, $D_F$ has no $L^2$-harmonic spinors,
for any face $F \subset M$, including $F = M$.
\end{corollary}

\begin{proof}\
The first part is an immediate consequence of Theorem
\ref{Theorem.c}. Since $0 \not \in \sigma_e(D_M)$, we have that $0
\in \sigma(D_M)$ if, and only if, $0$ is an eigenvalue of $D_M$.
\end{proof}

The following theorem takes care of the case when there are faces
of dimension zero, and hence completes our discussion.

\begin{theorem}\
Let \mate{M} is as in Theorem \ref{thm.LN} and $W \to M$ be a
Clifford bundle over $A^*$. Assume $M$ has faces of dimension
zero. Let $D$ be the Dirac operator with coefficients in $W$. Then
\begin{equation*}
    \sigma_e(D) = \RR.
\end{equation*}
\end{theorem}

\begin{proof}\ Use the same reasoning as in \cite{LN1}. Let $F$ be
a face of $M$ of dimension zero (that is, $F$ consists of one
point). The restriction $D_F$ of $D$ to the (subgroupoid
corresponding to the) face $F$ is the Dirac operator over $F
\times \RR^{n}$ with coefficients in the pull-back of $W \vert_F$
to $F \times \RR^{n}$, where $n$ is the dimension of $M$. Since
the spectrum of the Dirac operator on $\RR^{n}$ is $\RR$ (this can
be proved using the argument in the proof of Theorem
\ref{Theorem.c}), it follows from Theorem \ref{thm.Fredholm} (or
from \cite{LN1, MelrosePiazza}) that $\RR \subset \sigma_e(D)$.
\end{proof}

The results above extend to Dirac operators coupled with 
bounded potentials.

Our results on the spectrum of the Dirac operator are similar and
compatible with the results of \cite{BaerFV}, where the spectrum of
the Dirac operator on a manifold of finite volume is determined also
in terms of the properties of the boundary at infinity. The setting in
B\"ar's paper \cite{BaerFV} is different from ours (although
conformally equivalent). See also \cite{ammann.comin, BaerMHS,
BaerDS, MoroianuWL}.

\bibliographystyle{plain}
%\bibliography{prden}
\def\cprime{$'$} \def\cprime{$'$} \def\cprime{$'$} \def\cprime{$'$}
  \def\ocirc#1{\ifmmode\setbox0=\hbox{$#1$}\dimen0=\ht0 \advance\dimen0
  by1pt\rlap{\hbox to\wd0{\hss\raise\dimen0
  \hbox{\hskip.2em$\scriptscriptstyle\circ$}\hss}}#1\else {\accent"17 #1}\fi}

\end{document}